\documentclass[article]{amsart}

\usepackage{CJK,CJKnumb}
\usepackage{color}
\usepackage{indentfirst}
\usepackage{latexsym,bm}
\usepackage{amsmath,amssymb}
\usepackage{graphicx}
\usepackage{bbm}
\usepackage{amsthm}
\usepackage[ansinew]{inputenc}
\usepackage{array}
\usepackage{amsxtra}
\usepackage{amstext}
\usepackage{latexsym}
\usepackage{amsfonts}
\usepackage{bm}
\usepackage[pdftex,pdfstartview=FitH,bookmarksnumbered,bookmarksopen,bookmarksopenlevel=2,CJKbookmarks]{hyperref}

\newtheorem{theorem}{Theorem}[section]
\newtheorem{definition}[theorem]{Definition}

 \numberwithin{equation}{section}

\begin{document}

{\title[An integrable gradient flow for differential forms]{
{
An integrable example of gradient flows based on optimal transport of differential forms
}}
}
\author{
Yann Brenier and Xianglong Duan}
\date\today



\address{
CNRS UMR 7640 \\ 
Ecole Polytechnique   \\ 
Palaiseau\\
France}

\email{yann.brenier@polytechnique.edu and xianglong.duan@polytechnique.edu}

\begin{abstract}

Optimal transport theory has been a powerful tool for the analysis of parabolic equations
viewed as gradient flows of volume forms according to suitable transportation metrics.
In this paper, we present an example of gradient flows for closed $(d-1)-$forms in the
Euclidean space $\mathbb{R}^d$. In spite of its apparent complexity, the resulting very
degenerate parabolic system is fully integrable and can be viewed as the Eulerian
version of the heat equation for curves in the Euclidean space.
We analyze this system in terms of ``relative entropy" and ``dissipative solutions"
and provide global existence and weak-strong uniqueness results.
\end{abstract}
\maketitle

\section*{Introduction}

Optimal transport theory has been a powerful tool for the analysis of parabolic equations
viewed as gradient flows of volume forms according to suitable transportation metrics \cite{AGS,JKO,Ot,Sa,Vi}.
The theory of optimal transport for differential forms is not yet fully developed but there has been
some recent
progress, especially for symplectic forms and contact forms \cite{DGK,Re}. However, to the best
of our knowledge, little is known about
gradient flows in that context.
In this paper we just present an example of gradient flows for closed $(d-1)-$forms in the
Euclidean space $\mathbb{R}^d$. Such forms can be identified to divergence-free vector fields.
Our example, set on the flat torus $\mathbb{T}^d=\mathbb{(R/Z)}^d$, reads
\begin{equation}
\label{heat1}
\partial_t B+\nabla\cdot\left(\frac{B\otimes P-P\otimes B}{\rho}\right)=0,\;\;\;\nabla\cdot B=0,
\end{equation}
\begin{equation}
\label{heat2}
\partial_t\rho+\nabla\cdot P=0,\;\;\;
P=\nabla\cdot\left(\frac{B\otimes B}{\rho}\right)
\end{equation}
[or, in coordinates:
$$
\partial_t B^i+\partial_j\left(\frac{B^iP^j-P^iB^j}{\rho}\right)=0,\;\;\;\partial_j B^j=0,
$$
$$
\partial_t\rho+\partial_j P^j=0,\;\;\;
P^i=\partial_j\left(\frac{B^iB^j}{\rho}\right)],
$$
where $B$ is a time dependent divergence-free vector field (i.e. a closed $(d-1)-$form), $\rho$
is a time-dependent companion volume-form, and $P$ stands for $\rho v$, where $v$ is the time-dependent
velocity field transporting both $\rho$ and $B$ as differential forms.
As will be shown, this system turns out to be the gradient flow of functional
\begin{equation}
\label{GF1}
\mathcal{F}[\rho,B]=\int_{x\in \mathbb{T}^d} F(\rho(x),B(x)),\;\;\;
F(\rho,B)=\frac{|B|^2}{2\rho},
\end{equation}
according to the transportation metric
\begin{equation}
\label{GF2}
{||v||_\rho}=\sqrt{\int_{\mathbb{T}^d}v^2\rho},
\end{equation}
which is just the most usual transport metric for volume-forms \cite{AGS,Ot,Sa,Vi}.
Notice that $\mathcal{F}$, which is homogeneous of degree one, is a lower semi-continuous functional valued in $[0,+\infty]$ of
$\rho$ and $B$ viewed as Borel measures, respectively valued in $\mathbb{R}$
and $\mathbb{R}^d$ (cf. \cite{DT}), namely
$$
\mathcal{F}[\rho,B]=\sup\left\{
\int_{\mathbb{T}^d} \theta\rho+\Theta\cdot B,\;\;\;
\theta+\frac{1}{2}|\Theta|^2\le 0\right\}
\in [0,+\infty],
$$
where the supremum is taken over all $(\theta,\Theta)\in C(\mathbb{T}^d;\mathbb{R}\times\mathbb{R}^d)$.
\\
\\
From the PDE viewpoint, system (\ref{heat1},\ref{heat2}) is of degenerate parabolic type and can also be written,
in non-conservative form,
\begin{equation}
\label{heatnc}
\partial_t b+(v\cdot\nabla)b=(b\cdot\nabla)v,\;\;\;v=(b\cdot\nabla)b,
\end{equation}
or, equivalently,
\begin{equation}
\label{heatncbis}
\partial_t b=b\otimes b:D_x^2 b
\end{equation}
[In coordinates:
$$
\partial_t b^i+v^j\partial_j b^i=b^j\partial_j v^i,\;\;\;v^i=b^j\partial_j b^i
\;\;{\rm{or}}\;\;\partial_t b^i=b^k b^j\partial^2_{jk} b^i],
$$
where $b,v$ are the reduced variables $B/\rho,P/\rho$, and $\rho$ just solves the continuity equation
\begin{equation}
\partial_t\rho +\nabla\cdot(\rho v)=0.
\end{equation}
This system is formally integrable and can be viewed as the Eulerian
version of the heat equation for curves in the Euclidean space. More precisely, if
$(B,\rho,P)$ is of form
\begin{equation}
\label{smirnov}
(B,\rho,P)(t,x)=\int_{a\in\mathcal{A}}\left(\int_{\mathbb{R/Z}}\delta(x-X(t,s,a))
(\partial_s X,1,\partial_t X)(t,s,a)
ds\right) d\mu(a)
\end{equation}
where $(\mathcal{A},d\mu)$ is an abstract probability space of labels $a$, and, for $\mu-$a.e. label $a$, and every time $t$,
$s\in\mathbb{R/Z}\rightarrow X(t,s,a)\in \mathbb{R}^d$ is a loop in $\mathbb{T}^d$ (which means there is $N(a)\in\mathbb{Z}^d$
such that $X(t,s+1,a)=X(t,s,a)+N(a)$), subject to the heat equation
$$
\partial_t X(t,s,a)=\partial^2_{ss}X(t,s,a),
$$
then $(B,\rho,P)$ is expected to be a solution of (\ref{heat1},\ref{heat2}), as long as there
is no self- or mutual intersection of the different loops.
The goal of our paper is to provide a robust notion of generalized solutions which includes
fields $(B,\rho,P)$ of form (\ref{smirnov}) as global solutions of
system (\ref{heat1},\ref{heat2}), that we call, from now on, the Eulerian heat equation.

Our definition reads:
\begin{definition}\label{def:dddf1}
Let us fix $T>0$.
We say that $(B,\rho,P)$ with
$$(B,\rho)\in C([0,T],C(\mathbb{T}^d,\mathbb{R}^{d}\times\mathbb{R})'_{w^{*}}),\;\;\;P\in C([0,T]\times \mathbb{T}^d,\mathbb{R}^{d})',$$
is a dissipative solution of the Eulerian heat equation (\ref{heat1},\ref{heat2})
with initial data $(B_{0},\rho_0\ge0)\in
C(\mathbb{T}^d;\mathbb{R}^{d}
\times\mathbb{R})'
$
if and only if
\\
\\
{\rm\bf(i)} $\displaystyle{(B(0),\rho(0))=(B_{0},\rho_0)}$, $\;\;\nabla\cdot B=0,\;\;\partial_t\rho+\nabla\cdot P=0,$
in the sense of distributions;\\
{\rm\bf(ii)} For any $t\in[0,T]$, and any $b^{*},v^{*}\in C^{1}([0,T]\times\mathbb{T}^d,\mathbb{R}^d)$,
\begin{equation}
\label{eq:dddf-3}
E(t)e^{-rt}+\int^{t}_{0}e^{-rt'}\left[\int_{\mathbb{T}^d}  \frac{\widetilde{U}^{{\rm T}}Q_{r}\widetilde{U}}{2\rho}(t') - R(t')\right]{\rm d}t' \le
E(0)<+\infty.
\end{equation}
holds true for all $t\in [0,T]$ and $r\ge r_0$,
where
$$\widetilde{U}=\left(B-\rho b^{*},P-\rho v^{*}\right),
\;\;\;E=\int_{\mathbb{T}^d} \frac{{|B-\rho b^*|}^2}{2\rho},$$
$$Q_{r}=\left(
                       \begin{array}{cccc}
                          -\nabla v^{*}-\nabla {v^{*}}^{{\rm T}}+rI_{d}
&\nabla b^{*}-\nabla {b^{*}}^{{\rm T}}\\
                          \nabla {b^{*}}^{{\rm T}}-\nabla b^{*} & 2I_{d}\\
                       \end{array}
                     \right),$$
$$R=\int_{\mathbb{T}^d} \rho\;\mathrm{L}_{1}+ B\cdot\mathrm{L}_{2}+ P\cdot\mathrm{L}_{3},$$
$$\mathrm{L}_1={v^*}^2+D_t^*
\left(\frac{{b^*}^2}{2}\right)-(b^*\cdot\nabla)(b^*\cdot v^*),
\;\;\;D_t^*=(\partial_t+v^*\cdot\nabla),$$
$$\mathrm{L}_2=-D_t^* b^*+(b^*\cdot\nabla)v^*,
\;\;\;\mathrm{L}_3=-v^*+(b^*\cdot\nabla)b^*,$$
and $r_0$ is a constant depending on $b^*,v^*$, chosen so that $Q_{r_0}\geq I_{2d}$ for every $t,x$.

\end{definition}

This concept of solutions turns out to be convex in $(B,\rho,P)$ which is crucial
to include fields $(B,\rho,P)$ of form (\ref{smirnov}) as global dissipative solutions.
It is based on the relative entropy method, quite well known in the theory of
hyperbolic systems of conservation laws \cite{Da,KV,LT}, kinetic theory \cite{StR}, parabolic equations \cite{Ju}, and continuum mechanics
\cite{DST,FN}, just to quote few examples. It is closely related to
Lions' concept of dissipative solutions for the Euler equations of incompressible fluids \cite{Li} and related models
\cite{Brvla,BDS,Vo,Brtopo,BGL}. It is also related to the way the heat equation is solved for general metric measure spaces in \cite{AGS}.
With such a robust concept, we can see which way fields $(B,\rho,P)$
of form (\ref{smirnov}) are, indeed, solutions:

\begin{theorem}\label{theorem1}
Let
$$
\left(B_0,\rho_0\right)(x)=\int_{a\in\mathcal{A}}\left(\int_{\mathbb{R/Z}}\delta(x-X_0(s,a))\left(\partial_s X_0(s,a),1\right)ds\right) d\mu(a),\;\;\;x\in \mathbb{T}^d,
$$
where $(\mathcal{A},d\mu)$ is an abstract probability space of labels $a$, and
the loops
$s\in\mathbb{R/Z}\rightarrow X_0(s,a)\in \mathbb{T}^d$ are chosen so that
\begin{equation}
\label{cond}
\int_{x\in\mathbb{T}^d}
\frac{
|
\int_{\mathcal{A}\times\mathbb{R/Z}}
\delta(x-X(0,s,a))
\partial_s X(0,s,a)ds d\mu(a)|^2
}
{\int_{\mathcal{A}\times\mathbb{R/Z}}
\delta(x-X(0,s,a))ds d\mu(a)}
\end{equation}
$$
=\int_{\mathcal{A}\times\mathbb{R/Z}}
|\partial_s X(0,s,a)|^2ds d\mu(a)<+\infty.
$$
Then there is a global dissipative solution $(B,\rho,P)$ to the Eulerian heat equation (\ref{heat1},\ref{heat2}), explicitly given by
$$
(B,\rho,P)(t,x)=\int_{a\in\mathcal{A}}\left(\int_{\mathbb{R/Z}}\delta(x-X(t,s,a))
(\partial_s X,1,\partial_t X)(t,s,a)ds\right) d\mu(a),
$$
where, for $\mu-$a.e. $a$, $X(\cdot,\cdot,a)$ is solution to the heat equation
$$
\partial_t X(t,s,a)=\partial^2_{ss}X(t,s,a).
$$
\end{theorem}
Notice that condition (\ref{cond}) essentially means that, at time $t=0$, in the definition of
$(B_0,\rho_0)$, the loops $s\rightarrow X(0,s,a)$
have been chosen without self- or mutual intersections.
At this stage, we don't know how this result can be extended to $all$
Borel measures $B_0$, $\rho_0$ respectively valued in $\mathbb{R}^d$ and $\mathbb{R}_+$
such that
$$\nabla\cdot B_0=0,\;\;\;\int_{\mathbb{T}^d} \frac{|B_0|^2}{\rho_0}<+\infty.
$$
This is a delicate question of geometric measure theory, closely related to the topics discussed in \cite{AM,Sm}.
\\
\\
We also get the following ``weak-strong" uniqueness result:

\begin{theorem}\label{theorem2}

Suppose $(B,\rho,P)
\in C^{1}([0,T]\times\mathbb{T}^d,\mathbb{R}^d\times\mathbb{R}_{+}\times \mathbb{R}^d)$ is a classical solution to the Eulerian heat equation (\ref{heat1},\ref{heat2}).
Then $(B,\rho,P)$ is the unique dissipative solution to (\ref{heat1},\ref{heat2}) in the
sense of definition \ref{def:dddf1} with initial condition $(B_0,\rho_0)$.
\end{theorem}
 Here, the weak-strong uniqueness essentially means that any weak solution must coincide with a strong solution emanating from the same initial data as long as the latter exists. In other words, the strong solutions must be unique within the class of weak solutions. This kind of property has been widely studied in various kinds of PDEs (Navier-Stokes, Euler, etc.), even for measure valued solutions \cite{BDS}.

\subsection*{Miscellaneous remarks}

\subsubsection*{Remark 1}
By reducing a complicated degenerate parabolic system to the one-dimensional heat equation for loops, we follow, in a somewhat
opposite direction,  the path of Evans, Gangbo, Savin \cite{EGS}
who treated a degenerate parabolic system
in $\mathbb{R}^d$
with polyconvex entropy as an integrable system, by reducing it, in its
Eulerian version, to the scalar
equation in $\mathbb{R}^d$.

\subsubsection*{Remark 2}
In a companion paper \cite{BD}, the authors study the more delicate system
\begin{equation}\label{curve-shortening}
\partial_t B+\nabla\cdot\left(\frac{B\otimes P-P\otimes B}{|B|}\right)=0,\;\;\;\nabla\cdot B=0,
\;\;\;
P=\nabla\cdot\left(\frac{B\otimes B}{|B|}\right),
\end{equation}
where the algebraic constraint $\rho=|B|$ substitutes for the continuity equation.
Then, the analysis gets substantially more difficult. However, there is still some underlying integrability where the geometric heat equation
for curves (or curve-shortening flow, which corresponds to co-dimension $d-1$ mean curvature motion) substitutes for the linear
heat equation.
\subsection*{Acknowledgments}
This work has been partly supported by the ANR contract ISOTACE.
The first author is grateful to the Erwin Schr\"odinger Institute for its hospitality
during the first stage of this work. He also thanks the INRIA team MOKAPLAN where the work was
partly completed.

\section{Gradient flows for closed $(d-1)$-forms 	and transportation metrics}

Optimal transport theory is largely about giving a Riemannian structure to the space of volume forms
$\rho dx^1\wedge\cdot\cdot\cdot\wedge dx^d$
in $\mathbb{R}^d$
and many gradient flows can be derived accordingly, following the seminal work of Otto and collaborators \cite{JKO,Ot,Sa,Vi}.
Here we want to extend this theory to the case of  closed $d-1$ differential forms
in $\mathbb{R}^d$, or, in other words,  divergence-free vector fields. For instance, as $d=3$,
$$
B=B^1 dx^2\wedge dx^3+B^2 dx^3\wedge dx^1+B^3 dx^1\wedge dx^2,
$$
$$
0=dB=(\partial_1 B^1+\partial_2 B^2+\partial_3 B^3)dx^1\wedge dx^2\wedge dx^3={\rm div\;} B\;dx^1\wedge dx^2\wedge dx^3,
$$
and these formulae easily extend to arbitrary dimensions $d$.
For simplicity, we will only discuss about $\mathbb{Z}^d-$periodic forms so that we will use the flat torus
$\mathbb{T}^d=\mathbb{(R/Z)}^d$ instead of the entire space $\mathbb{R}^d$.

\subsection{Elementary closed $(d-1)-$forms and superposition of loops}
An elementary example of closed $(d-1)-$forms is given by
\begin{equation}
\label{form}
B(x)=\int_{\mathbb{R/Z}}\delta(x-X(s))X'(s)ds,\;\;\;x\in \mathbb{T}^d,
\end{equation}
where $X$
is a loop, i.e. a Lipschitz map $s\in \mathbb{R/Z}\rightarrow X(s)\in\mathbb{T}^d$
and $B$ should be understood (with an abuse of
notation) just as the vector-valued
distribution defined by
$$
<B^i,\beta_i>=
\int_{\mathbb{R/Z}}\beta(X(s))\cdot X'(s)ds,
$$
for all trial function $\beta\in C^\infty(\mathbb{T}^d;\mathbb{R}^d)$. In particular, $\nabla\cdot B=0$
immediately follows since
$$
<\partial_iB^i,\phi>=-<B^i,\partial_i\phi>
=-\int_{\mathbb{R/Z}}\nabla\phi(X(s))\cdot X'(s)ds
$$
$$
=-\int_{\mathbb{R/Z}}\frac{d}{ds}(\phi(X(s)))ds=0,
$$
for all trial function $\phi\in C^\infty(\mathbb{T}^ d;\mathbb{R})$. A much larger class of
divergence-free vector fields $B$ can be obtained by superposing loops:
\begin{equation}
\label{smirnov0}
B(x)=\int_{a\in\mathcal{A}}\left(\int_{\mathbb{R/Z}}\delta(x-X(s,a))
\partial_s X(s,a)ds\right) d\mu(a)
\end{equation}
where $(\mathcal{A},d\mu)$ is an abstract probability space of labels $a$, and, for $\mu-$a.e. label $a$,
$s\in\mathbb{R/Z}\rightarrow X(s,a)\in \mathbb{T}^d$ is a loop.
It is an important issue of geometric measure theory to see how large is this class
\cite{AM,Sm}. [A typical expected result being that
every divergence-free vector field $B$, with bounded mass
$$
\int_ {\mathbb{T}^d} |B|
=\sup \{<B,\beta>,\;\;\;\beta\in C^\infty(\mathbb{T}^d;\mathbb{R}^d),\;\;\;
\sup|\beta(x)|\le 1\;\}<+\infty,
$$
can be approximated by a superposition of $N$ loops $(X_\alpha,\alpha\in\{1,\cdot\cdot\cdot,N\})$ in the sense
$$
<B^i,\beta_i>=\lim_{N\rightarrow\infty}
\frac{1}{N}
\sum_{\alpha=1}^N
\int_{\mathbb{R/Z}}\beta(X_\alpha(s))\cdot X_\alpha'(s)ds,\;\;\;\forall
\beta\in C^\infty_0(\mathbb{T}^d;\mathbb{R}^d),
$$
$$
\int_ {\mathbb{T}^d} |B|=\lim_{N\rightarrow\infty}\frac{1}{N}\sum_{\alpha=1}^N
\int_{\mathbb{R/Z}}|X_\alpha'(s)|ds.]
$$
\subsection{Transportation of closed $(d-1)-$forms}

As for volume forms, the concept of transport involves time-dependent forms $B(t,x)$
and velocity fields $v(t,x)\in \mathbb{R}^d$. To get the transport equation right, we can refer to the case of
a moving loop $(t,x)\rightarrow X(t,s)\in \mathbb{T}^d$ subject to
$$
\partial_t X(t,s)=v(t,X(t,s)).
$$
Given a smooth trial function $\beta\in C^\infty(\mathbb{T}^d;\mathbb{R}^d)$, we get for the form attached to $X$
$$
B(t,x)=\int_{\mathbb{R/Z}}\delta(x-X(t,s))\partial_s X(t,s)ds,
$$
$$
\frac{d}{dt}<B^i(t,\cdot),\beta_i>=
\frac{d}{dt}\int_{\mathbb{R/Z}}\beta_i(X(t,s))\partial_s X^i(t,s)ds
$$
$$
=\int_{\mathbb{R/Z}}(\partial_j\beta_i)(X(t,s))\partial_t X^j(t,s)\partial_s X^i(t,s)ds
+\int_{\mathbb{R/Z}}\beta_i(X(t,s))\partial^2_{ts} X^i(t,s)ds
$$
$$
=\int_{\mathbb{R/Z}}(\partial_j\beta_i)(X(t,s))\partial_t X^j(t,s)\partial_s X^i(t,s)ds
-\int_{\mathbb{R/Z}}(\partial_j\beta_i)(X(t,s))\partial_s X^j(t,s)\partial_t X^i(t,s)ds
$$
$$
=\int_{\mathbb{R/Z}}(\partial_j\beta_i)(X(t,s))\left(v^j(t,X(t,s))\partial_s X^i(t,s)-\partial_s X^j(t,s)v^i(t,X(t,s))\right)ds
$$
$$
=<B^i v^j-B^j v^i,\partial_j\beta_i>
=-<\partial_j(B^i v^j-B^j v^i),\beta_i>
$$
So, for the transport of $B$ by $v$, we have found equation
$$
\partial_t B+\nabla\cdot(B\otimes v-v\otimes B)=0,
$$
which is linear in $B$ and stays valid, as a matter of fact, for all closed $(d-1)$-differential forms
(see \cite{CDK} for example).
In the case $d=3$, this equation is
usually called induction equation in the framework of (ideal) Magnetohydrodynamics,
with $B$ being interpreted as a magnetic field and $v$ as the velocity field of a conductive fluid. We will retain the name of
induction equation for any dimension $d$.
For the sequel of our discussion, it is convenient to attach to each time-dependent form $B$ generated by some loop $X$, a
companion volume form defined as
$$
\rho(t,x)=\int_{\mathbb{R/Z}}\delta(x-X(t,s))ds\ge 0,\;\;\;x\in \mathbb{T}^d,
$$
or, equivalently, by duality,
$$
<\rho(t,\cdot),\phi>=\int_{\mathbb{R/Z}}\phi(X(t,s))ds,\;\;\;\forall\phi\in C^\infty_0(\mathbb{T}^d;\mathbb{R}).
$$
[Notice that, in contrast with $B$, the definition of $\rho$ is affected by a change of parameterization of the loop $X$
with respect to $s$.]
The transportation equation for $\rho$ is nothing but
$$
\partial_t\rho+\nabla\cdot(\rho v)=0,
$$
just as in regular optimal transportation theory.
This equation is usually called continuity equation in Fluid Mechanics and we will also retain this name for any dimension $d$.
\subsection*{Transportation cost}

Mimicking the case of volume forms, which corresponds to regular optimal transportation,
we define a transportation cost by introducing,
for each fixed form $B$,  a Hilbert norm, possibly depending on $B$,
for all suitable transporting velocity field $x\in\mathbb{T}^d\rightarrow v(x)\in\mathbb{R}^d$.
In the case when we attach a volume form $\rho$ to $B$, the Hilbert norm may depend on both $B$ and $\rho$ and we denote
it by
$$
v\rightarrow ||v||_{B,\rho}.
$$
For volume forms, the most popular choice of norm is $v\rightarrow \sqrt{\int_{\mathbb{T}^d} |v|^2\rho}$
and we will concentrate on this choice in the present paper. Then, the resulting norm depends only on $\rho$ and is
simply denoted by $||\cdot||_\rho$.

\subsection*{Steepest descent}

Let us give a functional $\mathcal{F}[\rho,B]$ and monitor its steepest descent according to the Hilbert norm
$v\rightarrow ||v||_{B,\rho}$
on the space of velocity fields $v$ transporting $B$ and $\rho$. We will concentrate soon on the special case
when
$$
\mathcal{F}[\rho,B]=\int_{x\in\mathbb{T}^d}\frac{|B(x)|^2}{2\rho(x)},\;\;\;{||v||_\rho}
=\sqrt{\int_{x\in\mathbb{T}^d}|v(x)|^2\rho(x)dx},
$$
which turns out to be, in some sense, the simplest choice, as will be seen later on.
Nevertheless, let us start our calculations in the larger framework when
$$
\mathcal{F}[\rho,B]=\int_{x\in\mathbb{T}^d} F(\rho(x),B(x)),
$$
for some function $F:\mathbb{R}_+\times\mathbb{R}^d\rightarrow \mathbb{R}$, supposed to be
smooth away from $\rho=0$, such as $F(\rho,B)=\frac{B^2}{2\rho}$, for instance.
Thanks to the continuity and the induction equations, we get (assuming $B$, $\rho$ and $v$
to be smooth with $\rho>0$, using coordinates and denoting
by $F_\rho$ and $F_B$ the partial derivatives of $F$ with respect to $\rho$ and $B$)
$$
\frac{d}{dt}\mathcal{F}[\rho,B]=\int_{\mathbb{T}^d} \left(F_\rho \partial_t \rho+F_{B^i}\partial_t  B^i\right)
$$
$$
=-\int_{\mathbb{T}^d} \left[F_\rho \partial_i(\rho v^i)+F_{B^i}\partial_j (B^i v^j-B^j v^i)\right]
$$
(using the induction and continuity equations)
$$
=\int_{\mathbb{T}^d} \left[\rho\partial_i(F_\rho)-\big(\partial_j(F_{B^i})
-\partial_i(F_{B^j})\big)B^j\right]v^i
=
-\int_{\mathbb{T}^d} v\cdot\ G
$$
where
$$
G_i=-\rho\partial_i(F_\rho)+\big(\partial_j(F_{B^i})
-\partial_i(F_{B^j})\big)B^j.
$$
So, we have obtained
$$
\frac{d}{dt}\mathcal{F}[\rho,B]=-\int_{\mathbb{T}^d} v\cdot\ G
\ge  -\frac{1}{2}||v||^2_{\rho,B}-\frac{1}{2}{||G||^*_{\rho,B}}^2,
$$
where  ${||\cdot||^*}_{\rho,B}$ is the dual norm defined by
$$
\frac{1}{2}{||g||^*_{\rho,B}}^2=\sup_w\int_{\mathbb{T}^d} g\cdot w-\frac{1}{2}||w||^2_{\rho,B}.
$$
To get the steepest descent according to norm $||\cdot||_{\rho,B}$
it is enough to saturate this inequality so that
$$
\frac{1}{2}{||G||^*_{\rho,B}}^2+\frac{1}{2}||v||^2_{\rho,B}=\int G\cdot v,
$$
or, in other words, to define  $v$ as
the derivative with respect to $G$ of half the dual norm squared:
$$
v^i=\frac{d}{dG_i}\left(\frac{1}{2}{||G||^*_{\rho,B}}^2\right).
$$
\subsection*{The Eulerian version of the heat equation}
From now on, we will concentrate
on the special case (\ref{GF1},\ref{GF2}), namely
$$
\mathcal{F}[\rho,B]=\int_{x\in\mathbb{T}^d} F(\rho(x),B(x)),\;\;\;
F(\rho,B)=\frac{|B|^2}{2\rho},\;\;\;{||v||_\rho}=\sqrt{\int_{\mathbb{T}^d}v^2\rho}.
$$
According to the previous calculations, we first find
$$
F_\rho=-\frac{|B|^2}{2\rho^2},\;\;\;F_{B^i}=\frac{B_i}{\rho},
$$
$$
G_i=-\rho\partial_i(F_\rho)+\big(\partial_j(F_{B^i})
-\partial_i(F_{B^j})\big)B^j
$$
$$
=\rho\partial_i\left(\frac{B_j}{\rho}\right)\frac{B^j}{\rho}+\partial_j\left(\frac{B^i}{\rho}\right)B^j-\partial_i\left(\frac{{B_j}}{\rho}\right)B^j
=\partial_j\left(\frac{B^i B^j}{\rho}\right)
$$
(using $\partial_j B^j=0$).
Next, we get
$$
\frac{1}{2}{||G||^*_{\rho,B}}^2=\int\frac{G^2}{2\rho},
$$
and its derivative with respect to $G$ is just
${G}/{\rho}.$
Finally, we have obtained the steepest descent of $\mathcal{F}$
with respect to the transportation metric $v\rightarrow ||v||_\rho$,
precisely when
$$
v^i=\frac{1}{\rho}\partial_j(\frac{B^i B^j}{\rho}),
$$
which, combined with the induction and the continuity equations,
leads to the system (\ref{heat1},\ref{heat2}), namely
$$
\partial_t B+\nabla\cdot(\frac{B\otimes P-P\otimes B}{\rho})=0,\;\;\;\nabla\cdot B=0,
$$
$$
\partial_t\rho+\nabla\cdot P=0,\;\;\;
P=\nabla\cdot\left(\frac{B\otimes B}{\rho}\right)\;,
$$
where $P$ stands for $\rho v$ (i.e. the momentum, in terms of Fluid Mechanics).
When the solution $(B,\rho)$ is smooth with $\rho>0$
(which is definitely not the case when $B$ is generated by a single loop),
this system can be written in ``non conservative" form, in terms of the rescaled field $b=B/\rho$.
Indeed, from (\ref{heat1}), we get (in coordinates)
$$
0=b^i(\partial_t \rho+\partial_j(\rho v^j))+\rho\left(\partial_t b^i+v^j\partial_j b^i-b^j\partial_j v^i\right)
-v^i\partial_j(\rho b^j)
$$
$$
=\rho\left(\partial_t b^i+v^j\partial_j b^i-b^j\partial_j v^i\right)
$$
and
$$
v^i=b^j\partial_j b^i,
$$
which leads to the non-conservative version (\ref{heatnc})
of (\ref{heat1},\ref{heat2}), namely
$$
\partial_t b+(v\cdot\nabla)b=(b\cdot\nabla)v,\;\;\;v=(b\cdot\nabla)b,
$$
in which $\rho$ plays no role. This equation can also be written as the (very) degenerate parabolic system  (\ref{heatncbis}),
namely
$$
\partial_t b=b\otimes b:D_x^2 b.
$$
[Indeed, we get from (\ref{heatnc}), in coordinates,
$$
\partial_t b^i=-v^j\partial_j b^i+b^j\partial_j v^i
=-b^k\partial_k b^j\partial_j b^i+b^j\partial_j (b^k\partial_k b^i)
$$
$$
=-b^k\partial_k b^j\partial_j b^i+b^jb^k\partial^2_{kj} b^i+b^j\partial_j b^k\partial_{k} b^i=
b^jb^k\partial^2_{kj} b^i.]
$$

In a completely different direction, we can interpret (\ref{heat1},\ref{heat2}) just as a hidden
version of the one-dimensional heat equation!
Indeed, let us assume that a time dependent
loop $(t,s)\rightarrow X(t,s)$ solves the
linear heat equation
\begin{equation}
\label{heat}
\partial_t X(t,s)=\partial^2_{ss}X(t,s).
\end{equation}
and never self-intersects during some time interval $[0,T]$, so that we may find two smooth time-dependent vector field $v$
and $b$ such that
$$\partial_t X(t,s)=v(t,X(t,s)),\;\;\;\partial_s X(t,s)=b(t,X(t,s)).$$
Using the chain-rule, we first recover $v=(b\cdot\nabla)b$ directly from (\ref{heat}) and
then we observe that $\partial_t b+(v\cdot\nabla)b=(b\cdot\nabla)v$ is just the compatibility condition for $b$ and $v$ to be partial derivatives
of $X$. Surprisingly enough, we $directly$ recover the $non-conservative$ version of (\ref{heat1},\ref{heat2}).
We may recover the conservative form (\ref{heat1},\ref{heat2}) by reversing the computation we did to get the
non-conservative form, after $adding$ the field $\rho$ as solution of
the continuity equation  $\partial_t\rho+\nabla\cdot(\rho v)=0$, with initial condition
$$
\rho(0,x)=\int_{\mathbb{R/Z}}\delta(x-X(0,s))ds.
$$
Indeed, this equation is linear in $\rho$ and admits, since $v$ is supposed to be smooth, a unique distributional
solution which must be
$$
\rho(t,x)=\int_{\mathbb{R/Z}}\delta(x-X(t,s))ds
$$
since  $\partial_t X(t,s)=v(t,X(t,s))$.
To conclude this subsection, let us rename system (\ref{heat1},\ref{heat2}) as the Eulerian heat equation.

\subsubsection*{Remark 1}
At this stage, it seems strange to solve a complicated set of non-linear PDEs such as (\ref{heat1},\ref{heat2})
while we may, instead, solve the trivial one-dimensional heat equation! However, the derivation of (\ref{heat1},\ref{heat2})
we just performed is crucially based on the assumptions we made that $X$ is smooth (which is not a problem since
$X$ solves the heat equation) and non self-intersecting which, except in some very special situations, is not true
globally in time. So we can view (\ref{heat1},\ref{heat2}) as a way of extending the evolution beyond the first self-intersection time.
As a matter of fact, a similar situation is very well known for a collection of particles, labelled by some parameter
$a$,
and solving the trivial equation
$$
\partial^2_{tt} X(t,a)=0.
$$
Assuming the existence of a smooth velocity field $v$ such that
$$\partial_t X(t,a)=v(t,X(t,a)),$$
we immediately
obtain $\partial_t v+(v\cdot\nabla)v=0$ which is nothing but the multi-dimensional version of
the so-called inviscid Burgers equation, or, in other words, the non-conservative version of
the pressure-less Euler equations
\begin{equation}
\label{pressureless}
\partial_t P+\nabla\cdot\left(\frac{P\otimes P}{\rho}\right)=0,
\;\;\;
\partial_t\rho+\nabla\cdot P=0.
\end{equation}
As well known, in this model, particles of different labels may cross as time goes on (especially in the
case when parameter $a$ is continuously distributed). This is why system (\ref{pressureless}) is far from being
well understood, except in some special case, typically as $d=1$, where we may use the order of the real line to get a satisfactory formulation (such as in \cite{BGSW}).
Let us finally mention, as already done in the introduction, the work by Evans, Gangbo and Savin \cite{EGS} where the authors
are able to solve a complicated degenerate parabolic system with ``polyconvex entropy" in
$\mathbb{R}^d$ by noticing that its Eulerian
version is nothing but the regular scalar heat equation $\mathbb{R}^d$.

\subsubsection*{Remark 2}
We can also obtain the gradient flow structure for the Eulerian version of the curve-shortening flow (\ref{curve-shortening}). Forget about the volume form $\rho$, we only consider the transportation of the closed $(d-1)-$forms $B$ by vector fields $v$. Then (\ref{curve-shortening}) turns out to be the gradient flow of the following functional and transportation metric
$$\mathcal{F}[B]=\int_{x\in\mathbb{T}^d}|B(x)|,\;\;\;\|v\|_{B}=\sqrt{\int_{\mathbb{T}^d}v^2|B|}.$$


\section{Dissipative solutions to the Eulerian heat equation}
\label{main}
Let us consider a loop $X$ solution to the one-dimensional heat equation (\ref{heat})
and introduce the relative entropy
$$
\mathcal{E}(t)
=\int_{\mathbb{R/Z}}\frac{{|\partial_s X(t,s)-b^*(t,X(t,s))|}^2}{2}ds
$$
where $b^*\in C^\infty([0,T]\times\mathbb{T}^d;\mathbb{R}^d)$ is a fixed trial function.
We find, after elementary but lengthy computations, (see Appendix 1 for more details)
\begin{equation}\label{eq:in-1}
\begin{array}{r@{}l}
& \displaystyle{\quad
\frac{d\mathcal{E}}{dt}
=-\int_{\mathbb{R/Z}}|\partial_t X-v^*(t,X)|^2 ds
} \\
& \displaystyle{\quad
+\frac{1}{2}
\int_{\mathbb{R/Z}}(\partial_s X^i-{b^*}^i(t,X))(\partial_s X^j-{b^*}^j(t,X))
(\partial_j v^*_{i}+\partial_i v^*_{j})(t,X)ds
} \\
& \displaystyle{\quad \;\;\;\;
-\int_{\mathbb{R/Z}}
\big(\partial_s X^i-{b^*}^i(t,X)\big)\big(\partial_t X^j-{v^*}^j(t,X)\big)
(\partial_j b^*_{i}-\partial_i b^*_{j})(t,X)ds }\\
& \displaystyle{\quad
\;\;\;\;\;\;\;\;\;\;\;\;\;\;
+ \int_{\mathbb{R/Z}} \big( \mathrm{L}_1(t,X)+\partial_s X\cdot \mathrm{L}_2(t,X)+
\partial_t X\cdot \mathrm{L}_3(t,X)\big)ds
},\\
\end{array}
\end{equation}
where
$$\mathrm{L}_1={v^*}^2+D_t^*
\left(\frac{{b^*}^2}{2}\right)-(b^*\cdot\nabla)(b^*\cdot v^*),
\;\;\;D_t^*=(\partial_t+v^*\cdot\nabla),$$
$$\mathrm{L}_2=-D_t^* b^*+(b^*\cdot\nabla)v^*,\;\;\;\mathrm{L}_3=-v^*+(b^*\cdot\nabla)b^*.$$


In order to get more compact notations, we introduce
$$\widetilde{W}(t,s)=\big(\partial_s X(t,s)-b^{*}(t,X(t,s)),\partial_t X(t,s)-v^{*}(t,X(t,s))\big),$$
$$Q(b^{*},v^{*})=\left(
                       \begin{array}{cccc}
                          -\nabla v^{*}-\nabla {v^{*}}^{{\rm T}}
&\nabla b^{*}-\nabla {b^{*}}^{{\rm T}}\\
                          \nabla {b^{*}}^{{\rm T}}-\nabla b^{*} & 2I_{d}\\
                       \end{array}
                     \right)$$

Then \eqref{eq:in-1} can be written as
\begin{equation}\label{eq:in-2}
\frac{d\mathcal{E}}{dt}+ \int_{\mathbb{R/Z}}
 \frac{\widetilde{W}^{{\rm T}}Q(b^{*},v^{*})\widetilde{W}}{2}(t,s)ds -\mathcal{R}(t)=0,
\end{equation}
where
$$
\mathcal{R}(t)=\int_{\mathbb{R/Z}} \big( \mathrm{L}_1(t,X(t,s))
+\partial_s X(t,s)\cdot \mathrm{L}_2(t,X(t,s))+ \partial_t X(t,s)\cdot \mathrm{L}_3(t,X(t,s))\big)ds
$$

We use $I_{n:m}$ to represent the $n\times n$ diagonal matrix whose first $m$ terms are 1 while the rest terms are 0, let $I_{d}$ be the $d\times d$ identity matrix. Then we can choose $r_0\geq0$ big enough, in terms
of the trial functions $b^*,v^*$, such that $$Q_{r_0}=Q(b^{*},v^{*})+r_0I_{2d:d}\geq I_{2d}>0.$$
In addition, we observe that
$$
\frac{1}{2}{(\widetilde{W}^T I_{2d:d}\;\widetilde{W})(t,s)}
=\frac{1}{2}{\big|\partial_s X(t,s)-b^*(t,X(t,s))\big|}^2,
$$
which is exactly the relative entropy density.
Thus, for any constant $r\geq r_0$, we obtain by integrating in time (\ref{eq:in-2})
after multiplication by $e^{-rt}$:
\begin{equation}\label{eq:in-3}
\mathcal{E}(t)e^{-rt}+\int^{t}_{0}e^{-rt'}\left[\int_{\mathbb{R/Z}}
\frac{1}{2}(\widetilde{W}^{{\rm T}}Q_{r}\widetilde{W})(t',s)ds - \mathcal{R}(t')\right]{\rm d}t'=\mathcal{E}(0)
\end{equation}
Let us now consider a collection of loops, labelled by $a\in\mathcal{A}$, where $(\mathcal{A},d\mu)$ is an abstract probability space,
and subject to the heat equation
$$
\partial_t X(t,s,a)=\partial^2_{ss}X(t,s,a).
$$
We set, for each $a\in\mathcal{A}$,
$$
\widetilde{W}(t,s,a)=\big(\partial_s X(t,s,a)-b^{*}(t,X(t,s,a)),\partial_t X(t,s,a)-v^{*}(t,X(t,s,a))\big),
$$
$$
\mathcal{R}(t,a)=\int_{\mathbb{R/Z}} \big( \mathrm{L}_1(t,X)
+\partial_s X\cdot \mathrm{L}_2(t,X)+ \partial_t X\cdot \mathrm{L}_3(t,X)\big)(t,s,a)ds,
$$
$$
\mathcal{E}(t,a)
=\int_{\mathbb{R/Z}}\frac{{|\partial_s X(t,s,a)-b^*(t,X(t,s,a))|}^2}{2}ds,
$$
so that (\ref{eq:in-3}) reads
\begin{equation}\label{eq:in-3bis}
\mathcal{E}(t,a)e^{-rt}+\int^{t}_{0}e^{-rt'}\left[\int_{\mathbb{R/Z}}
\frac{1}{2}(\widetilde{W}^{{\rm T}}Q_{r}\widetilde{W})(t',s,a)ds - \mathcal{R}(t',a)\right]{\rm d}t' =\mathcal{E}(0,a).
\end{equation}
Next, we introduce the averaged fields $(B,\rho,P)$:
$$
(B,\rho,P)(t,x)=\int_{a\in\mathcal{A}}\left(\int_{\mathbb{R/Z}}\delta(x-X(t,s,a))
(\partial_s X,1,\partial_t X)(t,s,a)ds\right) d\mu(a)
$$
and, also,
$$
\widetilde{U}(t,x)=\big(B(t,x)-\rho(t,x) b^{*}(t,x),P(t,x)-\rho(t,x) v^{*}(t,x)\big)
$$
$$
=\int_{a\in\mathcal{A}}\left(\int_{\mathbb{R/Z}}\delta(x-X(t,s,a))
\big(\partial_s X-b^*(t,X),\partial_t X-v^*(t,X)\big)(t,s,a)ds\right) d\mu(a),
$$
$$R=\int_{\mathbb{T}^d}   \rho\;\mathrm{L}_{1}+ B\cdot\mathrm{L}_{2}
+ P\cdot\mathrm{L}_{3}
=\int_{a\in\mathcal{A}}\mathcal{R}(\cdot,a)d\mu(a),
\;\;\;
E=\int_{\mathbb{T}^d}\frac{{|B-\rho b^*|}^2}{2\rho}.
$$
We see that
$$
E(t)
=\int_{\mathbb{T}^d}\frac{{B}^2}{2\rho}(t)+
\int_{a\in\mathcal{A}}\int_{\mathbb{R/Z}}
\left(
-\partial_s X\cdot b^*(t,X)+\frac{1}{2}|b^*(t,X)|^2
\right)(t,s,a)
ds d\mu(a),
$$
$$
=\int_{a\in\mathcal{A}}\mathcal{E}(t,a)d\mu(a)
+\int_{\mathbb{T}^d}\frac{{B}^2}{2\rho}(t)
-\frac{1}{2}\int_{a\in\mathcal{A}}\int_{\mathbb{R/Z}}
|\partial_s X(t,s,a)|^2ds d\mu(a).
$$
By the Cauchy-Schwarz inequality
$$
\int_{\mathbb{T}^d}\frac{{B}^2}{\rho}(t)
=\int_{x\in \mathbb{T}^d}
\frac{
|\int_{a\in\mathcal{A}}\int_{\mathbb{R/Z}}\delta(x-X(t,s,a))
\partial_s X(t,s,a)ds d\mu(a)|^2
}
{\int_{a\in\mathcal{A}}
\int_{\mathbb{R/Z}}
\delta(x-X(t,s,a))ds d\mu(a)}
$$
$$
\le \int_{a\in\mathcal{A}}\int_{\mathbb{R/Z}}
|\partial_s X(t,s,a)|^2ds d\mu(a),
$$
so that
$$
E(t)\le \int_{a\in\mathcal{A}}\mathcal{E}(t,a)d\mu(a).
$$
In a similar way,
$$
\int_{\mathbb{T}^d}  \frac{\widetilde{U}^{{\rm T}}Q_{r}\widetilde{U}}{\rho}(t)
\le
\int_{a\in\mathcal{A}}\int_{\mathbb{R/Z}}
(\widetilde{W}^{{\rm T}}Q_{r}\widetilde{W})(t,s,a)ds d\mu(a).
$$
Thus, when integrating $equality$ (\ref{eq:in-3bis}) in $a\in\mathcal{A}$ with respect to $\mu$,
we deduce the following $inequality$
$$
E(t)e^{-rt}+\int^{t}_{0}e^{-rt'}\left[
\int_{\mathbb{T}^d}  \frac{\widetilde{U}^{{\rm T}}Q_{r}\widetilde{U}}{2\rho}(t')
- R(t')\right]dt' \le E(0),
$$
provided
$$
E(0)=\int_{a\in\mathcal{A}}\mathcal{E}(0,a)d\mu(a),
$$
which means that
the Cauchy-Schwarz inequality we used saturates at time $t=0$, i.e.
$$
\int_{x\in\mathbb{T}^d}
\frac{
|\int_{a\in\mathcal{A}}\int_{\mathbb{R/Z}}\delta(x-X(0,s,a))
\partial_s X(0,s,a)ds d\mu(a)|^2
}
{\int_{a\in\mathcal{A}}
\int_{\mathbb{R/Z}}
\delta(x-X(0,s,a))ds d\mu(a)}
$$
$$
=\int_{a\in\mathcal{A}}
|\partial_s X(0,s,a)|^2ds d\mu(a).
$$
This essentially means, as already explained in the Introduction,
that, at time $t=0$, in the definition of the initial fields
$(B,\rho)(t=0,\cdot)$, the loops $s\rightarrow X(0,s,a)$
have been chosen without self- or mutual intersections.
The resulting convex $inequality$ is precisely the one we have chosen
in the Introduction to define dissipative solutions
for the Eulerian heat equation (\ref{heat1},\ref{heat2}), through Definition \ref{def:dddf1}.

\subsection{Proof of Theorem \ref{theorem1}}

The proof has just been provided, while obtaining the concept of dissipative solutions!

\subsection{Proof of Theorem \ref{theorem2}}

The proof is straightforward. Suppose $(B,\rho,P)$ is a classical solution. Then it is easy to verify that the non-conservative variables $b=B/\rho,v=P/\rho$ solve 
$$
D_t b=(b\cdot\nabla)v,
\;\;\;D_t=(\partial_t+v\cdot\nabla),
\;\;\;v=(b\cdot\nabla)b,
$$
from which we easily deduce
$$
v^2+D_t
\left(\frac{v^2}{2}\right)=(b\cdot\nabla)(b\cdot v).
$$
Then, it is enough to set $b^*=b$ and $v^*=v$ in definition \ref{def:dddf1},
to make $L_1=L_2=L_3=0$, $R=0$ and $E(0)=0$.
So for any dissipative solution  $(B',\rho',P')$
in the sense of definition \ref{def:dddf1}, by the inequality (\ref{eq:dddf-3}), we have that
$$E(t)=0,\;\;\int^{T}_{0}\int_{\mathbb{T}^d}  e^{-r_0t}\frac{\widetilde{U}^{{\rm T}}Q_{r_0}\widetilde{U}}{2\rho}=0,$$
So we must have $B'=\rho' b$ and $P'=\rho' v$ since $Q_{r_0}\geq I_{2d}$. Now since both $\rho$ and $\rho'$ solve the same continuity equation
$$\partial_t\rho'+\nabla\cdot (\rho'v)=0.$$
with the same initial data $\rho'(0)=\rho(0)=\rho_0$, we mush have that $\rho'=\rho$, which completes the proof.


\section{Appendix 1: direct recovery of equation (\ref{eq:in-1})}

Let loop $X(t,s)$ be a solution to the heat equation (\ref{heat}). For any smooth vector field $b^*$, the relative entropy
$$
\mathcal{E}(t)
=\int_{\mathbb{R/Z}}\frac{{|\partial_s X(t,s)-b^*(t,X(t,s))|}^2}{2}ds$$
can be written as
$$
\mathcal{E}(t)=\int_{\mathbb{R/Z}}\left(\frac{|\partial_s X|^2}{2} -\partial_s X\cdot b^*(t,X) +\frac{|b^*(t,X)|^2}{2}\right) ds=\mathcal{E}_1 (t)+\mathcal{E}_2 (t)+\mathcal{E}_3 (t).
$$
For $\mathcal{E}_1 (t)=\int |\partial_s X|^2/2$, we have
$$\frac{d}{dt}\mathcal{E}_1(t)=\int \partial_s X\cdot\partial_{ts} X =-\int \partial_{ss}X \cdot \partial_t X=-\int|\partial_t X|^2.$$
(since $\partial_t X=\partial_{ss} X$). For any smooth vector field $v^*$, we have
$$-\int|\partial_t X|^2= \int -|\partial_t X-v^*(t,X)|^2 + |v^*(t,X)|^2-2\partial_t X\cdot v^*(t,X)$$
$$=\int -|\partial_t X-v^*(t,X)|^2 + |v^*(t,X)|^2-\partial_t X\cdot v^*(t,X)-\partial_{ss} X\cdot v^*(t,X).$$
In coordinates, we have
$$-\int\partial_{ss} X^i v^*_i(t,X)= \int \partial_s X^i \partial_{j}v^*_i(t,X)\partial_s X^j$$
$$=\frac{1}{2}\int\big(\partial_s X^i-{b^*}^i(t,X)\big)\big(\partial_s X^j-{b^*}^j(t,X)\big)(\partial_{j}v^*_i + \partial_{i}v^*_j)(t,X)$$
$$+\int \partial_sX^i{b^*}^j(t,X)(\partial_{j}v^*_i + \partial_{i}v^*_j)(t,X) - \big({b^*}^i{b^*}^j\partial_jv^*_i\big)(t,X). $$
So we have,
$$\frac{d}{dt}\mathcal{E}_1(t)=\int-|\partial_t X-v^*(t,X)|^2 + \mathrm{L}'_1(t,X)+\partial_s X\cdot \mathrm{L}'_2(t,X)+
\partial_t X\cdot \mathrm{L}'_3(t,X) $$
$$+\frac{1}{2}\int\big(\partial_s X^i-{b^*}^i(t,X)\big)\big(\partial_s X^j-{b^*}^j(t,X)\big)(\partial_{j}v^*_i + \partial_{i}v^*_j)(t,X),$$
where
$$\mathrm{L}'_1=|v^*|^2-{b^*}^i{b^*}^j\partial_jv^*_i,\;\;\;(\mathrm{L}'_2)_i=(\partial_{j}v^*_i + \partial_{i}v^*_j){b^*}^j,\;\;\; (\mathrm{L}'_3)_i=-v^*_i .$$
Now let's look at $\mathcal{E}_2(t)=-\int\partial_s X \cdot b^*(t,X)$. In coordinates, we have,
$$\frac{d}{dt}\mathcal{E}_2(t)=\int-\partial_{st} X^i \cdot b^*_i(t,X)-\partial_s X^i(\partial_t b^*_i)(t,X)-\partial_s X^i\partial_t X^j (\partial_j b^*_i)(t,X)$$
$$=\int -\partial_s X^i\partial_t X^j(\partial_j b^*_i-\partial_i b^*_j)(t,X) -\partial_s X^i(\partial_t b^*_i)(t,X) $$
$$=-\int \big(\partial_s X^i-{b^*}^i(t,X)\big)\big(\partial_t X^j-{v^*}^j(t,X)\big)(\partial_j b^*_i-\partial_i b^*_j)(t,X)$$
$$+ \int\mathrm{L}''_1(t,X)+\partial_s X\cdot \mathrm{L}''_2(t,X)+
\partial_t X\cdot \mathrm{L}''_3(t,X),$$
where
$$\mathrm{L}''_1=(\partial_j b^*_i-\partial_i b^*_j){b^*}^i{v^*}^j,\;\;(\mathrm{L}''_2)_i=-\partial_t b^*_i-(\partial_j b^*_i-\partial_i b^*_j){v^*}^j,\;\;(\mathrm{L}''_3)_i=(\partial_j b^*_i-\partial_i b^*_j){b^*}^j.$$
For the last term $\mathcal{E}_3(t)=\int |b^*(t,X)|^2/2$, we have,
$$\frac{d}{dt}\mathcal{E}_3(t)= \int \big(\partial_t b^*_i(t,X) + \partial_t X^j\partial_j b^*_i(t,X)\big){b^*}^i(t,X).$$
So, in summary, we have
$$\frac{d\mathcal{E}}{dt}=\frac{1}{2}\int\big(\partial_s X^i-{b^*}^i(t,X)\big)\big(\partial_s X^j-{b^*}^j(t,X)\big)(\partial_{j}v^*_i + \partial_{i}v^*_j)(t,X)$$
$$-\int \big(\partial_s X^i-{b^*}^i(t,X)\big)\big(\partial_t X^j-{v^*}^j(t,X)\big)(\partial_j b^*_i-\partial_i b^*_j)(t,X)$$
$$\int-|\partial_t X-v^*(t,X)|^2 + \mathrm{L}_1(t,X)+\partial_s X\cdot \mathrm{L}_2(t,X)+
\partial_t X\cdot \mathrm{L}_3(t,X),$$
where
$$\mathrm{L}_1=\mathrm{L}'_1+\mathrm{L}''_1 + \partial_t\left(\frac{{b^*}^2}{2}\right)={v^*}^2+D_t^*
\left(\frac{{b^*}^2}{2}\right)-(b^*\cdot\nabla)(b^*\cdot v^*),
\;\;D_t^*=(\partial_t+v^*\cdot\nabla)$$
$$\mathrm{L}_2=\mathrm{L}'_2+\mathrm{L}''_2=-D_t^* b^*+(b^*\cdot\nabla)v^*+\nabla(b^*\cdot v^*)$$
$$\mathrm{L}_3=\mathrm{L}'_3+\mathrm{L}''_3 + \nabla\left(\frac{{b^*}^2}{2}\right)=-v^*+(b^*\cdot\nabla)b^*.$$
Since
$$\int\partial_s X\cdot \big[\nabla(b^*\cdot v^*)\big](t,X)=\int \partial_{s}\big(b^*(t,X)\cdot v^*(t,X)\big)=0,$$
we can remove the gradient term $\nabla(b^*\cdot v^*)$ from $\mathrm{L}_2$ and finally get (\ref{eq:in-1}).


\end{document}